# Alternating Euler sums at the negative integers

Khristo N. Boyadzhiev,
Ohio Northern University
Department of Mathematics
Ada, Ohio, 45810, USA
k-boyadzhiev@onu.edu

H. Gopalkrishna Gadiyar and R. Padma
AU-KBC Research Centre
M. I. T. Campus of Anna University
Chromepet, Chennai 600 044, INDIA
gadiyar@au-kbc.org
padma@au-kbc.org

This paper is dedicated to the memory of Professor R. Sitaramachandra Rao

**Abstract**. We study three special Dirichlet series, two of them alternating, related to the Riemann zeta function. These series are shown to have extensions to the entire complex plane and we find their values at the negative integers (or residues at poles). These values are given in terms of Bernoulli and Euler numbers.

## 1. Introduction

Let

$$H_n = 1 + \frac{1}{2} + \frac{1}{3} + \ldots + \frac{1}{n}, \tag{1.1}$$

($n = 1, 2, \ldots$) be the harmonic numbers. Apostol-Vu [2] and Matsuoka [9] studied the function

$$h(s) = \sum_{n=1}^{\infty} \frac{1}{n^s} H_n, \tag{1.2}$$

and showed that it has a second-order pole at $s = 1$ and simple poles at zero and the negative odd integers. The values at the negative even integers were evaluated in terms of Bernoulli numbers.

These results were discussed and specified further in [6].

In this article we study three similar series

$$u(s) = \sum_{n=1}^{\infty} \frac{(-1)^{n-1}}{n^s} H_n, \qquad (1.3)$$

$$v(s) = \sum_{n=1}^{\infty} \frac{(-1)^{n-1}}{n^s} H_n^-, \qquad (1.4)$$

$$w(s) = \sum_{n=1}^{\infty} \frac{1}{n^s} H_n^-, \qquad (1.5)$$

where

$$H_n^- = 1 - \frac{1}{2} + \frac{1}{3} + \ldots + \frac{(-1)^{n-1}}{n}. \qquad (1.6)$$

The functions $u$ and $v$ are well defined and holomorphic for $\mathrm{Re}\,s > 0$, as $H_n \sim \log n$ when $n \to \infty$. The function $w$ is defined by (1.5) for $\mathrm{Re}\,s > 1$. We prove here that $u(s)$ extends to an entire function on the complex plane $\mathbb{C}$, $w(s)$ extends to $\mathbb{C}$ with a simple pole at $s = 1$, while $v(s)$ extends to $\mathbb{C}$ with poles at $0, -1, -3, \ldots$, much like $h(s)$. We also evaluate these functions (when defined) for $s = 0, -1, -2, \ldots$ The values $v(-2n)$ are connected to a certain convolution of Bernoulli numbers in the spirit of [6].

In our work we shall use the Riemann zeta function ($\mathrm{Re}\,s > 1$)

$$\zeta(s) = \sum_{n=1}^{\infty} \frac{1}{n^s}, \qquad (1.7)$$

and the Euler-Dirichlet eta function ($\mathrm{Re}\,s > 0$)

$$\eta(s) = \sum_{n=1}^{\infty} \frac{(-1)^{n-1}}{n^s}, \qquad (1.8)$$

with the relationship between them for all complex numbers

$$\eta(s) = (1 - 2^{1-s})\zeta(s). \qquad (1.9)$$

The alternating Euler sums above have been studied and evaluated for certain positive integers by Sitaramachandra Rao [11] in terms of zeta values. See also [8].

**2. The Euler-Boole summation formula.**

In order to show that the function $h(s)$ extends to all complex numbers and to investigate its values at the negative integers, Apostol-Vu and Matsuoka used the Euler-Maclaurin summation formula. For the alternating series $u(s)$ and $v(s)$ it is appropriate to use the alternating version of this formula, called in [3] the Euler-Boole formula (see also [10, p. 34]). The alternating version is less known, but equally important. While the ordinary Euler summation formula is based on the Bernoulli numbers and polynomials, the alternating version is naturally based on the Euler polynomials $E_n(x)$ defined by the generating function [1],

$$\frac{2e^{xt}}{e^t+1} = \sum_{n=0}^{\infty} E_n(x)\frac{t^n}{n!} \quad (|t|<\pi). \tag{2.1}$$

Thus

$$E_0(x)=1,\ E_1(x)=x-\frac{1}{2},\ E_2(x)=x^2-x,\ E_3(x)=x^3-\frac{3}{2}x^2+\frac{1}{4},\ \text{etc}. \tag{2.2}$$

The following lemma can be found in [3].

**Lemma 1.** Let $f$ be a function defined for $t \geq x$ with continuous derivatives of order $q$ and such that $f^{(k)}(t) \to 0\ (t\to\infty), k=0,1,\ldots,q$. Then

$$\sum_{k=0}^{\infty}(-1)^k f(x+k) = \frac{1}{2}\sum_{m=0}^{q-1}\frac{E_m(0)}{m!}f^{(m)}(x)+R_q, \tag{2.3}$$

where

$$R_q = \frac{1}{2}\int_0^{\infty}\overline{E}_{q-1}(-t)\frac{f^{(q)}(x+t)}{(q-1)!}dt,$$

and $\overline{E}_q$ is an antiperiodic function defined as follows

$$\overline{E}_q(x) = E_q(x)\ \text{for}\ 0\leq x\leq 1,\ \overline{E}_q(x+1) = -\overline{E}_q(x)\ \text{etc}, \tag{2.4}$$

( $\overline{E}_q$ is periodic with period 2 ).

Note that

$$E_n(0) = \frac{2}{n+1}(1 - 2^{n+1})B_{n+1},  \quad (2.5)$$

where $B_k$ are the Bernoulli numbers. Thus

$$E_{2n}(0) = 0, \; n = 1, 2, \ldots . \quad (2.6)$$

(The infinite version of (2.3) without remainder can be seen in [4, (3.30)].)

Taking $f(x) = x^{-s}$, $Re\, s > 0$, we find

$$f^{(m)}(x) = (-1)^m s(s+1)\ldots(s+m-1)\frac{1}{x^{s+m}}. \quad (2.7)$$

Set here $x = n$. Then (2.3) gives (with some $q > 1 + |Re\, s|$ )

$$\sum_{k=0}^{\infty} \frac{(-1)^k}{(n+k)^s} = \frac{1}{2n^s} - \frac{1}{2}\sum_{m=0}^{q-1} \frac{s(s+1)(s+2)\ldots(s+2m)}{n^{s+2m+1}(2m+1)!} E_{2m+1}(0)$$

$$+ \frac{s(s+1)(s+2)\ldots(s+2q-1)}{(2q-1)!} \int_0^{\infty} \frac{\bar{E}_{2q-1}(-t)dt}{2(n+t)^{s+2q}}. \quad (2.8)$$

**Theorem 1**. The function $u(s)$ admits an extension to the closed left half plane $Re\, s \leq 0$ as an entire function. Moreover,

$$u(0) = \frac{1}{2}\ln 2, \quad (2.9)$$

and for $n \geq 1$,

$$2u(-2n) = \eta(1-2n) + nE_{2n-1}(0), \quad (2.10)$$

or

$$u(-2n) = \frac{(2^{2n}-1)(1-2n)}{4n}B_{2n}.$$

Also,

$$u(-1) = \frac{1}{4} - \frac{1}{4}\ln 2 \qquad (2.11)$$

and for $n \geq 2$

$$2u(1-2n) = E_1(0)(2n-1)\eta(3-2n) \qquad (2.12)$$

$$+ \frac{1}{3!}E_3(0)(2n-1)(2n-2)(2n-3)\eta(5-2n)$$

$$+ \frac{1}{5!}E_5(0)(2n-1)(2n-2)(2n-3)(2n-4)(2n-5)\eta(7-2n)$$

$$+ \ldots + E_{(2n-1)}(0)\eta(1)$$

(In the last term here, $(2n-1)!$ appears in both the numerator and the denominator and thus cancels out.)

*Proof.* By changing the order of summation we have

$$u(s) = \sum_{n=1}^{\infty} \frac{(-1)^{n-1}}{n^s} \sum_{k=1}^{n} \frac{1}{k} = \sum_{n=1}^{\infty} \frac{1}{n} \sum_{k=n}^{\infty} \frac{(-1)^{k-1}}{k^s} \qquad (2.13)$$

$$= \sum_{n=1}^{\infty} \frac{(-1)^{n-1}}{n} \sum_{k=0}^{\infty} \frac{(-1)^k}{(n+k)^s}.$$

The last sum can be evaluated by using (2.8). Multiplying (2.8) by $\dfrac{(-1)^{n-1}}{n}$ and summing from $n=1$ to infinity we obtain

$$2u(s) = \eta(s+1) - \sum_{m=0}^{q-1} \frac{s(s+1)(s+2)\ldots(s+2m)}{(2m+1)!} E_{2m+1}(0)\eta(s+2m+2)$$

$$+ \frac{s(s+1)(s+2)\ldots(s+2q-1)}{(2q-1)!} \int_0^{\infty} \overline{E}_{2q-1}(-t)\phi^-(s+2q,t)\,dt, \qquad (2.14)$$

where
$$\phi^-(s,t) = \sum_{n=1}^{\infty} \frac{(-1)^{n-1}}{n(n+t)^s}. \tag{2.15}$$

This important representation provides a lot of information about the function $u(s)$. It shows that this function extends to $Re\ s \leq 0$. For every complex number $s$ we can take large enough $q$ to show that $u(s)$ becomes well defined. For the evaluation at the negative integers we remind that, according to (1.7)

$$\eta(-p) = (1 - 2^{1+p})\zeta(-p), p = 0, 1, \ldots \tag{2.16}$$

and since

$$\zeta(0) = \frac{-1}{2},\ \zeta(-2n) = 0,\ \zeta(1-2n) = -\frac{B_{2n}}{2n}, n = 1, 2, \ldots \tag{2.17}$$

we have

$$\eta(0) = \frac{1}{2},\ \eta(-2n) = 0,\ \eta(1-2n) = \frac{2^{2n}-1}{2n}B_{2n}, n = 1, 2, \ldots \tag{2.18}$$

Next, setting $s = 0, -1, -2$ we find

$$u(0) = \frac{1}{2}\eta(1) = \frac{1}{2}\ln 2, \tag{2.19}$$

$$2u(-1) = \eta(0) - E_1(0)(-1)\eta(1),\ \text{i.e.,}$$

$$u(-1) = \frac{1}{4}(1 - \ln 2), \tag{2.20}$$

$$2u(-2) = \eta(-1) - E_1(0)(-2)\eta(0),$$

$$u(-2) = \frac{-1}{8}. \tag{2.21}$$

Following this pattern, we obtain from (2.14) the general formula

$$2u(-2n) = \eta(1-2n) + 2n E_{2n-1}(0)\eta(0) \tag{2.22}$$

which gives (2.10).

The function $u(s)$ was evaluated for the positive even integers by Sitaramachandra Rao [11]. For instance,

$$u(2) = \frac{5}{8}\zeta(3). \tag{2.23}$$

**Remark 1.** Instead of the function $u(s)$ one can consider the more general function

$$u(s;p) = \sum_{n=1}^{\infty} \frac{(-1)^{n-1}}{n^s}(1 + \frac{1}{2^p} + \ldots + \frac{1}{n^p}) = \sum_{n=1}^{\infty} \frac{(-1)^{n-1}}{n^p} \sum_{k=0}^{\infty} \frac{(-1)^k}{(n+k)^s}, \tag{2.24}$$

where $p \geq 1$. The above proof shows that the only change in the representation (2.14) will be a shift $s \to s+p-1$ in the functions $\eta(s)$ and $\phi^-(s)$. All other formulas will undergo the same shift. Therefore, the essential case is $p = 1$. This remark applies also to the functions $v(s)$ and $w(s)$.

**Remark 2.** The numbers

$$G_n = 2(1 - 2^n)B_n = nE_{n-1}(0) \tag{2.25}$$

are known as the Genocchi [Jenoki] numbers. They are defined by the generating function

$$\frac{2t}{e^t + 1} = \sum_{n=0}^{\infty} G_n \frac{t^n}{n!} \quad (|t| < \pi) \tag{2.26}$$

## 3. The function $v(s)$

**Theorem 2.** The function $v(s)$ has the representation

$$2v(s) = \zeta(s+1) - \sum_{m=0}^{q-1} \frac{s(s+1)(s+2)\ldots(s+2m)}{(2m+1)!} E_{2m+1}(0) \zeta(s+2m+2)$$

$$+ \frac{s(s+1)(s+2)\ldots(s+2q-1)}{(2q-1)!} \int_0^{\infty} \bar{E}_{2q-1}(-t) \phi^+(s+2q, t)\, dt, \tag{3.1}$$

where $q > 1 + |\operatorname{Re} s|$ is a positive integer and

$$\phi^+(s,t) = \sum_{n=1}^{\infty} \frac{1}{n(n+t)^s}. \qquad (3.2)$$

Therefore, it admits a holomorphic extension to the closed left half plane $\operatorname{Re} s \le 0$ with simple poles at $s = 0, -1, -3, \ldots$. The residues at these poles are

$$\operatorname{Res}[s=0] = \frac{1}{2}, \qquad (3.3)$$

$$\operatorname{Res}[s=-2n-1] = \frac{1}{2} E_{2n+1}(0), \qquad (3.4)$$

$(n = 0, 1, \ldots)$.

Moreover, for $n = 1, 2, \ldots$,

$$2\nu(-2n) = \zeta(1-2n) - n E_{2n-1}(0). \qquad (3.5)$$

*Proof.* Changing the order of summation we obtain

$$\nu(s) = \sum_{n=1}^{\infty} \frac{(-1)^{n-1}}{n^s} \sum_{k=1}^{n} \frac{(-1)^{k-1}}{k} = \sum_{n=1}^{\infty} \frac{1}{n} \sum_{k=0}^{\infty} \frac{(-1)^k}{(n+k)^s}, \qquad (3.6)$$

and then the representation (3.1) follows from (2.8).

### 4. The function $w(s)$

**Theorem 3.** The function $w(s)$ defined by (1.6) for $\operatorname{Re} s > 1$ has the representation

$$w(s) = \frac{\eta(s)}{s-1} + \frac{1}{2}\eta(s+1) + \sum_{m=1}^{q} \frac{s(s+1)\ldots(s+2m-2)}{(2m)!} B_{2m}\eta(s+2m) + R_q(s), \qquad (4.1)$$

$$R_q(s) = -\frac{s(s+1)\ldots(s+2q)}{(2q+1)!} \int_0^{\infty} \bar{B}_{2q-1}(t)\, \phi^-(s+2q+1, t)\, dt, \qquad (4.2)$$

where $B_{2m}$ are the Bernoulli numbers, $q$ is a sufficiently large positive integer and $\bar{B}_k(t)$ are bounded functions. Thus $w(s)$ extends as an analytic function to all complex numbers with a simple pole at $s = 1$ with residue $\eta(1) = \log 2$. Its values $w(-p)$, $p = 0, 1, 2, \ldots$ can be found from the above representation (note that $R_q(-p) = 0$ when $2q > p$). In particular,

$$w(0) = \frac{1}{2}\eta(1) - \eta(0), \tag{4.3}$$

$$w(-1) = \frac{-1}{2}\eta(-1) + \frac{1}{2}\eta(0) - \frac{B_2}{2!}\eta(1), \tag{4.4}$$

$$w(-2) = \frac{1}{2}\eta(-1) - B_2\eta(0), \tag{4.5}$$

etc. The general formulas are

$$w(-2n) = \frac{1}{2}\eta(1-2n) - B_{2n}\eta(0), \tag{4.6}$$

$$w(1-2n) = \frac{-1}{2n}\eta(1-2n)$$

$$- \sum_{m=1}^{n} \frac{(2n-1)(2n-2)\ldots(2n-2m+1)}{(2m)!} B_{2m}\eta(-2n+2m+1), \tag{4.7}$$

for $n = 1, 2, \ldots$.

*Proof.* Changing the order of summation (just like in (2.13) and (3.5)) we write

$$w(s) = \sum_{n=1}^{\infty} \frac{1}{n^s} \sum_{k=1}^{n} \frac{(-1)^{k-1}}{k} = \sum_{n=1}^{\infty} \frac{(-1)^{n-1}}{n} \sum_{k=0}^{\infty} \frac{1}{(n+k)^s}. \tag{4.8}$$

The second sum on the right hand side is not alternating and we need to use the Euler-Maclaurin summation formula (as Apostol and Vu did on p.88 in [2]),

$$\sum_{k=0}^{\infty} \frac{1}{(n+k)^s} = \frac{1}{s-1}\frac{1}{n^{s-1}} + \frac{1}{2n^s} + \sum_{m=1}^{q} \frac{s(s+1)\ldots(s+2m-2)}{(2m)!} \frac{B_{2m}}{n^{s+2m-1}}$$

$$- \frac{s(s+1)\ldots(s+2q)}{(2q+1)!} \int_0^{\infty} \frac{\overline{B}_{2q-1}(t)dt}{(n+t)^{s+2q+1}}, \tag{4.9}$$

where $q$ is a large enough positive integer and $\bar{B}_k(t)$ is the periodic extension of the Bernoulli polynomial $B_k(t)$ from $[0,1]$ to $\mathbb{R}$ with period 1.

Multiplying (4.9) by $\dfrac{(-1)^{n-1}}{n}$ and summing for $n \geq 1$ we obtain the desired representation (4.1).

## 5. Equation for $v(s)$ with a Hankel contour integral

Our next theorem shows a connection between the function $v(s)$ defined by (1.4) and a special Hankel integral. For $s \neq 1, 2, \ldots$ define

$$G(s) = \frac{\Gamma(1-s)}{2\pi i} \int_L \frac{z^{s-1} e^z}{1 + e^z} \log\left(\frac{e^z - 1}{z}\right) dz, \qquad (5.1)$$

where the contour $L = L_- \cup L_+ \cup L_\epsilon$ consists of three parts: $L_-$ is the "lower side" (i.e. $arg(z) = -\pi$) of the ray $(-\infty, -\epsilon)$, $\epsilon > 0$, traced left to right, $L_+$ is the "upper side" ($arg(z) = \pi$) of this ray traced right to left, and $L_\epsilon = \{z = \epsilon e^{\theta i} : -\pi \leq \theta \leq \pi\}$ is a small circle traced counterclockwise and connecting the two sides of that ray. This contour was used also in [6].

The following theorem is a modification of the main result in [6].

**Theorem 4**. For every complex $s$ where both sides are defined we have the equation

$$G(s) = v(s) - \zeta(s+1) - \psi(s)\eta(s) - \eta'(s) \qquad (5.2)$$

(here $\psi(s)$ is the digamma function and the eta function is defined in (1.8)).

*Proof.* For $Re\, s > 1$ set

$$I(s) = \frac{1}{2\pi i} \int_L \frac{z^{s-1} e^z}{1 + e^z} \log\left(\frac{e^z - 1}{z}\right) dz. \qquad (5.3)$$

Now let $\epsilon \to 0$. The integral over $L_\epsilon$ turns into zero, as the function

$$\frac{e^z}{1+e^z} \log\left(\frac{e^z-1}{z}\right) \tag{5.4}$$

is holomorphic in a neighborhood of zero. Noticing that $z = xe^{-\pi i}$ on $L_-$ and $z = xe^{\pi i}$ on $L_+$, we find that

$$I(s) = \frac{e^{-\pi i s}}{2\pi i} \int_\infty^0 \frac{x^{s-1} e^{-x}}{1+e^{-x}} \log\left(\frac{1-e^{-x}}{x}\right) dx + \frac{e^{\pi i s}}{2\pi i} \int_0^\infty \frac{x^{s-1} e^{-x}}{1+e^{-x}} \log\left(\frac{1-e^{-x}}{x}\right) dx \tag{5.5}$$

$$= \frac{\sin \pi s}{\pi} \int_0^\infty \frac{x^{s-1} e^{-x}}{1+e^{-x}} \log\left(\frac{1-e^{-x}}{x}\right) dx. \tag{5.6}$$

At the same time,

$$\int_0^\infty \frac{x^{s-1} e^{-x}}{1+e^{-x}} \log\left(\frac{1-e^{-x}}{x}\right) dx = \int_0^\infty \frac{x^{s-1} e^{-x}}{1+e^{-x}} \log(1-e^{-x}) dx - \int_0^\infty \frac{x^{s-1}}{e^x+1} \log x\, dx, \tag{5.7}$$

and the two integrals on the right hand side are easy to evaluate. First, we use the representation

$$\frac{\log(1-e^{-x})}{1+e^{-x}} = \sum_{n=1}^\infty (-1)^n H_n^- e^{-nx}. \tag{5.8}$$

Multiplying both sides here by $x^{s-1} e^{-x}$ and integrating from zero to infinity we find

$$\int_0^\infty \frac{x^{s-1} e^{-x}}{1+e^{-x}} \log(1-e^{-x}) dx = \sum_{n=1}^\infty (-1)^n H_n^- \int_0^\infty x^{s-1} e^{-(n+1)x} dx$$

$$= \Gamma(s) \sum_{n=1}^\infty \frac{(-1)^n H_n^-}{(n+1)^s} = \Gamma(s)(v(s) - \zeta(s+1)). \tag{5.9}$$

Second, differentiating for $s$ the well-known representation

$$\Gamma(s)\eta(s) = \int_0^\infty \frac{x^{s-1}}{1+e^x} dx, \tag{5.10}$$

we obtain

$$\int_0^\infty \frac{x^{s-1}}{1+e^x} \log x \, dx = \Gamma'(s)\eta(s) + \Gamma(s)\eta'(s) = \Gamma(s)(\psi(s)\eta(s) + \eta'(s)). \tag{5.11}$$

Collecting together the two evaluations (5.9) and (5.11) we obtain

$$\int_0^\infty \frac{x^{s-1}e^{-x}}{1+e^{-x}} \log\left(\frac{1-e^{-x}}{x}\right) dx = \Gamma(s)(v(s) - \zeta(s+1) - \psi(s)\eta(s) - \eta'(s)), \tag{5.12}$$

i.e.
$$I(s) = \frac{1}{\pi}\Gamma(s)\sin(\pi s)(v(s) - \zeta(s+1) - \psi(s)\eta(s) - \eta'(s)), \tag{5.13}$$

and the theorem follows from here in view of the identity

$$\Gamma(s)\Gamma(1-s) = \frac{\pi}{\sin \pi s}. \tag{5.14}$$

**6. A special identity for Bernoulli numbers.**

In this section we shall find a relation between the values $v(-2n)$ and the numbers $C_n$, $n = 1, 2, \ldots$ defined as the convolutions.

$$C_n = \frac{1}{2}\sum_{k+j=n} \frac{1}{j!} E_j(1) \frac{B(k)}{k!k}, \quad k = 1, 2, \ldots; j = 0, 1, \ldots \tag{6.1}$$

where $B(n) = B_n$ are the Bernoulli numbers for $n \neq 1$ and $B(1) = -B_1 = 1/2$. Thus for $|z| < \pi$.

$$\frac{ze^z}{e^z - 1} = \frac{-z}{e^{-z} - 1} = \sum_{n=0}^\infty \frac{B(n)}{n!} z^n. \tag{6.2}$$

The notation $B(n)$ is used here in order to avoid the negative sign in $B_1$. Integrating (6.2) we find

$$\log\left(\frac{e^z - 1}{z}\right) = \sum_{n=1}^\infty \frac{B(n)}{n!n} z^n. \tag{6.3}$$

From (2.1)

$$\frac{e^z}{e^z+1} = \frac{1}{2}\sum_{n=0}^{\infty} \frac{1}{n!} E_n(1) z^n. \tag{6.4}$$

The product of the functions in (6.3) and (6.4) is the generating function of the numbers $C_n$,

$$\frac{e^z}{e^z+1} \log(\frac{e^z-1}{z}) = \sum_{n=1}^{\infty} C_n z^n. \tag{6.5}$$

The Euler polynomials have the property

$$E_n(t+1) + E_n(t) = 2t^n, \tag{6.6}$$

which shows that

$$E_n(1) = -E_n(0), \quad n \geq 1 \tag{6.7}$$

(and also $E_0(1) = E_0(0) = 1$ from (2.2)).

**Corollary 1.** For $n = 1, 2, \ldots$ we have

$$(2n)! C_{2n} = (\frac{1}{4n} + 2^{2n-1} - \frac{1}{2}) B_{2n} \tag{6.8}$$

*Proof.* This follows from Theorem 2 (equation (3.5)) and Theorem 4. From (5.2)

$$G(-2n) = v(-2n) - \zeta(1-2n) - (\psi(s)\eta(s) + \eta'(s))|_{s=-2n} \tag{6.9}$$

where the last term is, in fact, zero (cf [6, Corollary 1]). This is so because we can use the equation

$$\psi(s) = \psi(1-s) - \pi \cot \pi s \tag{6.10}$$

and write

$$\psi(s)\eta(s) = \psi(1-s)\eta(s) - \eta(s)\pi \cot \pi s. \tag{6.11}$$

For $s = -2n$ we have $\psi(1+2n)\eta(-2n) = 0$ and

$$\eta(s)\pi \cot \pi s|_{s=-2n} = \eta'(-2n). \tag{6.12}$$

as follows from the Taylor expansion around $s = -2n$,

$$\eta(s)\pi\cot\pi s = \eta'(-2n) + \frac{1}{2}\eta''(-2n)(s+2n) + O((s+2n)^2). \tag{6.13}$$

Therefore,

$$G(-2n) = v(-2n) - \zeta(1-2n). \tag{6.14}$$

At the same time, the numbers $G(-2n)$ are proportional to the Taylor coefficients of the series (6.5), i.e.

$$G(-2n) = \frac{(2n)!}{2\pi i} \int_{L_\epsilon} \frac{e^z}{1+e^z} \log\left(\frac{e^z-1}{z}\right) \frac{dz}{z^{2n+1}} = (2n)!\, C_{2n}. \tag{6.15}$$

This representation results from (5.1). For this special integrand the integrals on $L_+$ and $L_-$ cancel out and (6.16) follows from Cauchy's formula. From this we find immediately that

$$(2n)!\, C_{2n} = v(-2n) - \zeta(1-2n), \tag{6.16}$$

and (6.8) follows in view of (2.5), (2.17) and (3.5).